\documentclass[12pt]{scrartcl}
\usepackage{amsthm, amsmath, amssymb}
\usepackage[alphabetic]{amsrefs}
\usepackage{aurical}
\usepackage[T1,OT2,OT1]{fontenc}
\newcommand\cyr{%
\renewcommand\rmdefault{wncyr}%
\renewcommand\sfdefault{wncyss}%
\renewcommand\encodingdefault{OT2}%
\normalfont
\selectfont}
\DeclareTextFontCommand{\textcyr}{\cyr}

\usepackage{graphicx}
\newtheorem{theorem}{Theorem}
\newtheorem{corollary}[theorem]{Corollary}

\newtheorem{conj}[theorem]{Conjecture}

\newcommand{\R}{\mathbb{R}}

\newcommand{\e}{\varepsilon}

\renewcommand{\epsilon}{\varepsilon}

\renewcommand{\leq}{\leqslant}
\renewcommand{\le}{\leqslant}
\renewcommand{\geq}{\geqslant}
\renewcommand{\ge}{\geqslant}

\title{Long-range \\
phase coexistence models:\\
recent progress\\ on the fractional \\Allen-Cahn equation
\thanks{It is a pleasure to thank Xavier Cabr\'e and Joaquim Serra
for their very useful comments and for several pleasant discussions.
This work has been supported by the Australian
Research Council Discovery Project ``N.E.W. Nonlocal Equations at Work''.
The authors are members of GNAMPA/INdAM.}}

\author{Serena Dipierro and
Enrico Valdinoci}

\date{} 

\begin{document}

\maketitle 

\begin{abstract}
In this set of notes, we present some recent developments
on the fractional Allen-Cahn equation
$$ (-\Delta)^s u = u-u^3,$$
with special attention to $\Gamma$-convergence results,
energy and density estimates, convergence of level sets, Hamiltonian estimates,
rigidity and symmetry results.
\end{abstract}


The study of nonlocal models
for phase coexistence equations
is an interesting, and remarkably challenging, research topic 
which has experienced a rapid growth in the recent mathematical literature.
The goal of this paper is to collect several results
and present them in a unified and easily accessible way, with a style which
tries to combine, as much as possible, rigorous presentations and
intuitive descriptions of the problems under consideration and of
the methods used in some of the proofs.

The model will be also somehow described ``from scratch'' and, in spite of
the necessary simplifications which make the topic mathematically treatable,
we hope that we managed to preserve some important treats
from the physical model in view of the applications in material sciences,
using also these real-world motivations as a hint towards the development of
rigorous and quantitative theories.

In what follows, we will discuss specifically
$\Gamma$-convergence results,
energy and density estimates, convergence of level sets, Hamiltonian estimates,
rigidity and symmetry results.
To this end, we will first recall the classical Allen-Cahn phase coexistence model
in Section~\ref{S1}, where we will also present some of the classical results
and conjectures about it. Then, we will consider the nonlocal counterpart
of these problems in Section~\ref{S2}.

\section{The classical Allen-Cahn equation}\label{S1}

The so-called Allen-Cahn equation is
a semilinear, scalar equation, originally introduced by
John W. Cahn and Sam Allen in the 1970s.
In the stationary case, this equation is of elliptic type and can be written
in the form
\begin{equation}\label{AC1}
-\Delta u = u-u^3\quad{\mbox{ in }}\Omega.
\end{equation}
In~\eqref{AC1}, the function~$u=u(x)$
represents an
order parameter that determines the phase of the medium at a given point~$x\in\Omega$.
In this setting, the ``pure phases'' are denoted by the state parameters~$-1$
and~$1$, and the Laplacian term can be considered
as a surface tension or interfacial energies, which ends up preventing abrupt phase
changes from point to point and ``wild'' phase oscillations.
The set~$\Omega\subseteq\R^n$ can be viewed as the ``container''
and then equation~\eqref{AC1} aims at giving a simple, but effective,
description of phase coexistence.

One sees that equation~\eqref{AC1} possesses a variational structure
and solutions of~\eqref{AC1}
correspond to the critical points of the energy functional
\begin{equation}\label{E1}\begin{gathered}
{{ {\hbox{{\Fontamici J}}} }}_\Omega(u):=
\int_\Omega \left( \frac12\,|\nabla u(x)|^2
+W(u(x))\right)\,dx,\\
{\mbox{where }}\;
W(r):=\frac{(1-r^2)^2}{4}.\end{gathered}
\end{equation}
It is of course tempting to look at the ``big picture'' offered by this scenario:
namely, given a candidate function~$u_\star$,
one can consider the blow-down sequence
\begin{equation}\label{BLSDOW}
u_\e(x):=u_\star\left( \frac{x}\e\right)\qquad{\mbox{as }}\e\searrow0,
\end{equation}
and one remarks that
\begin{equation*}\begin{split} 
{{ {\hbox{{\Fontamici J}}} }}_{\Omega/\e}(u_\star)\,&=
\int_{\Omega/\e} \left( \frac{\e^2}2 \,|\nabla u_\e(\e x)|^2
+W\big( u_\e(\e x)\big) \right)\,dx\\&=\frac1{\e^n}\,
\int_{\Omega} \left( \frac{\e^2}{2}\,|\nabla u_\e(x)|^2
+W\big( u_\e (x)\big)\right)\,dx
\\&=\frac1{\e^{n-1}}\, {{ {\hbox{{\Fontamici J}}} }}_{\Omega,\e}(u_\e),\end{split}
\end{equation*}
where
\[ {{ {\hbox{{\Fontamici J}}} }}_{\Omega,\e}(v):= 
\int_{\Omega} \left( \frac{\e}2 \,|\nabla v(x)|^2
+\frac1\e\,W\big(v (x)\big)\right)\,dx.\]
In particular, $u_\star$ is a local
minimizer for~${{ {\hbox{{\Fontamici J}}} }}_{\Omega/\e}$ 
(that is, ${{ {\hbox{{\Fontamici J}}} }}_{\Omega/\e}(u_\star)
\le{{ {\hbox{{\Fontamici J}}} }}_{\Omega/\e}(u_\star+\varphi)$
for all~$\varphi\in C^\infty_0(\Omega/\e)$)
if and only if~$u_\e$
is a local minimizer for~${{ {\hbox{{\Fontamici J}}} }}_{\Omega,\e}$
(that is, ${{ {\hbox{{\Fontamici J}}} }}_{\Omega,\e}(u_\e)\le
{{ {\hbox{{\Fontamici J}}} }}_{\Omega,\e}(u_\e+\varphi)$
for all~$\varphi\in C^\infty_0(\Omega)$).
In this setting, the following results are classical:

\begin{itemize}
\item 
{{\emph{$\Gamma$-convergence.}}}
If~$\Omega$ is a smooth domain and~$u_\e:\Omega\to[-1,1]$ is a
sequence of minimizers for~${{ {\hbox{{\Fontamici J}}} }}_{\Omega,\e}$ such that
$$\sup_{\e\in(0,1)}{{ {\hbox{{\Fontamici J}}} }}_{\Omega,\e}(u_\e)<+\infty,$$
then, up to a subsequence,
\begin{equation}\label{limit E} \lim_{\e\searrow0} u_\e=\chi_E-
\chi_{\R^n\setminus E}\quad{\mbox{ in }}
L^1(\Omega),\end{equation}
for some set~$E\subseteq\R^n$ which minimizes the perimeter in~$\overline\Omega$
with respect to its boundary datum, see~\cite{MR0445362}.
\item {{\emph{Energy and density estimates.}}}
If~$R\ge1$ and~$
u:B_{R+1}\to[-1,1]$ is a minimizer of~${{ {\hbox{{\Fontamici J}}} }}_{B_{R+1}}$,
then
\begin{equation}\label{EN:DE1}
{{ {\hbox{{\Fontamici J}}} }}_{B_{R+1}}(u)\le CR^{n-1},
\end{equation}
for some~$C>0$.

In addition, if~$u(0)=0$, then
\begin{equation}\label{EN:DE2}\begin{gathered}
{\mbox{the Lebesgue measures of $\{u>1/2\}$ and~$\{u<-1/2\}$
in~$B_R$}}\\{\mbox{are both greater than~$cR^n$,}}\end{gathered}
\end{equation}
for some~$c>0$,
see~\cite{MR1310848}.
\item {{\emph{Locally uniform convergence of level sets.}}}
If~$\Omega$ is a smooth domain, $E\subseteq\R^n$ and~$u_\e:\Omega\to[-1,1]$
is a minimizer
of~${{ {\hbox{{\Fontamici J}}} }}_{\Omega,\e}$ such that~\eqref{limit E} holds true,
then the set~$\{ |u_\e|\le1/2\}$ converges locally uniformly in~$\Omega$
to~$\partial E$ as~$\e\searrow0$,
namely
$$ \lim_{\e\searrow0} \,\sup_{x\in\Omega'} {\operatorname{dist}}(x,\partial E) =0$$
for any~$\Omega'\Subset\Omega$,
see again~\cite{MR1310848}.
\item {{\emph{Pointwise gradient bounds.}}}
If~$u:\R^n\to[-1,1]$ is a critical point for the energy functional
in~\eqref{E1} for any bounded domain~$\Omega\subset\R^n$, then we have the 
following pointwise gradient
bound:
\begin{equation}\label{MODICA}
|\nabla u(x)|^2 \le 2W\big( u(x)\big)\quad{\mbox{for all }}x\in\R^n,
\end{equation}
see~\cite{MR803255}.
\end{itemize}

The inequality in~\eqref{MODICA}
can be seen as part of a family of formulas related
to Hamiltonian Identities, see~\cite{MR2381198}.
As a matter of fact,
we observe that in dimension~$1$, the bound in~\eqref{MODICA}
reduces to the classical Conservation of Energy Principle:
indeed, if~$u:\R\to[-1,1]$ is a solution of~$\ddot u=W'(u)$ in~$\R$,
it follows that
$$ \frac{d}{dx} \left(\frac{|\dot u(x)|^2}{2}-W(u(x))\right)=
\dot u(x)\,\ddot u(x)-
W'(u(x))\,\dot u(x)=0,$$
and therefore
\begin{equation}\label{823d91} \frac{|\dot u(x)|^2}{2}-W(u(x))
=\frac{|\dot u(y)|^2}{2}-W(u(y))\le \frac{|\dot u(y)|^2}{2},\end{equation}
for any~$y\in\R$.
Now, two cases occur: either~$\dot u$ never vanishes, or
\begin{equation}\label{CA:MODICA}
\dot u(y_0)=0\end{equation}
for some~$y_0\in\R$. In the first case, $u$ is monotone,
say increasing, and since it is bounded it has a limit at~$+\infty$
and, as a consequence,
$$ \lim_{y\to+\infty} \dot u(y)=0.$$
Using this information in~\eqref{823d91}, one obtains that
$$ \frac{|\dot u(x)|^2}{2}-W(u(x))\leq\lim_{y\to+\infty}\frac{|\dot u(y)|^2}{2}=0,$$
which is~\eqref{MODICA} in this case. If instead~\eqref{CA:MODICA} holds
true, it is enough to compute~\eqref{823d91} at~$y:=y_0$
and deduce~\eqref{MODICA} in this case as well.\medskip

Furthermore, one of the most important problems related to equation~\eqref{AC1}
is the following celebrated conjecture by Ennio De Giorgi:

\begin{conj}[see~\cite{MR533166}]\label{GDD}
Let~$u:\R^n\to[-1,1]$ be a solution of~\eqref{AC1}
in the whole of~$\R^n$ such that
\begin{equation}\label{MONP}
\frac{\partial u}{\partial x_n}(x)>0 \quad{\mbox{ for all }}x\in\R^n.
\end{equation}
Is it true that~$u$ is $1$D -- that is, $u(x)=u_0(\omega\cdot x)$ for some~$u_0:\R\to\R$
and~$\omega\in S^{n-1}$ -- at least if~$n\le8$?
\end{conj}

For a survey on Conjecture~\ref{GDD}, we refer to~\cite{MR2528756}.
See also~\cite{MR2569331} for a series of related rigidity
problems in elliptic equations. Here, we just recall that
a positive answer to Conjecture~\ref{GDD} (and, in fact, to more
general problems) has been provided in dimensions~$2$ and~$3$,
see~\cites{MR1637919, MR1655510, MR1775735, MR1843784}.

In its full generality and to the best of our knowledge,
the problem posed in Conjecture~\ref{GDD} is still open
in dimensions~$4$ to~$8$,
see~\cite{MR1954269}, but the claim in Conjecture~\ref{GDD} holds true also
in dimensions~$4$ to~$8$ under the following limit assumption:
\begin{equation}\label{LIMINU} \lim_{x_n\to\pm\infty} u(x',x_n)=\pm1,\end{equation}
see~\cites{MR2480601, MR2757359, MR3695374}.\medskip

Furthermore, condition~\eqref{LIMINU} can be relaxed.
For instance, in~\cite{MR2728579} the behavior of the limit profiles
$$ u^\pm(x'):=\lim_{x_n\to\pm\infty} u(x',x_n)$$
for monotone solutions~$u$ is taken into account, and it is proved
that if both~$ u^+$ and~$ u^-$ are~$ 2$D (i.e. they depend on at most two
Euclidean variables), then~\eqref{LIMINU} is automatically satisfied
and, as a consequence, if~$N\le8$, then~$ u$ is~$ 1$D.
Interestingly, when~$N\le4$, it is enough to suppose
that one between~$u^+$ and~$u^-$ is~$2$D in order to deduce~\eqref{LIMINU}
and that~$u$ is~$1$D (see in particular Theorems~1.1 and~1.2
in~\cite{MR2728579}).\medskip

These types of results focus on the $1$D symmetry ``inherited from the profiles~$u^\pm$ at infinity'' of the monotone solutions of~\eqref{AC1}.
Other symmetry results can be obtained from suitable asymptotic
properties of flatness type
of the level sets of the solutions, or from the graph properties
of such level sets, see in particular~\cites{MR2728579, MR3488250}.
Results for quasi-minimal solutions are also available in the literature,
see~\cites{MR2413100, MR2728579}.\medskip

Conjecture~\ref{GDD} is also related to a number of problems arising
in related fields, such as geometry and dynamical systems, see~\cite{MR2955315}
and the references therein.
\medskip

We also recall that a variant of Conjecture~\ref{GDD}, in which the monotonicity
property in~\eqref{MONP}
is replaced by a uniform limit assumption at infinity, was proposed
independently by
Gary William Gibbons:

\begin{conj}[see~\cite{MR1340266}] \label{GDD2}
Let~$u:\R^n\to[-1,1]$ be a solution of~\eqref{AC1}
in the whole of~$\R^n$ such that
\begin{equation}\label{LIMIU} \lim_{x_n\to\pm\infty} u(x',x_n)=\pm1\quad
{\mbox{ uniformly with respect to $x':=(x_1,\dots,x_{n-1})\in\R^{n-1}$}}. \end{equation}
Is it true that~$u(x)=u_0(x_n)$ for some~$u_0:\R\to\R$?
\end{conj}

Interestingly,
Conjecture~\ref{GDD2} turns out to be true in any dimension, as proved
independently in~\cites{MR1765681, MR1763653, MR1755949}. In this
sense, the uniform limit condition in~\eqref{LIMIU}
happens to be significantly stronger and to impose further crucial
restrictions when compared with the monotonicity assumption in~\eqref{MONP}.
Also, the uniform assumption in~\eqref{LIMIU}
ends up being dramatically stronger than the simple limit condition in~\eqref{LIMINU}
and so it provides significantly different types of results,
in particularly ruling out the example in~\cite{MR2846486}.\medskip

A variational variant of Conjecture~\ref{GDD} replaces
the monotonicity assumption in~\eqref{MONP} with a minimality assumption:
with respect to this, it is proved in~\cites{MR2480601}
that
\begin{equation}\label{89019230234855757757583}
\begin{gathered}
{\mbox{if $u:\R^n\to[-1,1]$ is a minimizer
for~${{ {\hbox{{\Fontamici J}}} }}_{\Omega}$}}\\{\mbox{in any bounded domain~$\Omega\subset\R^n$, and~$n\le7$, then $u$ is $1$D.
}}\end{gathered}\end{equation}
The example in~\cite{MR2846486} implies that a similar result cannot hold
in dimension~9 and higher. The case of dimension~8 has been recently
addressed in~\cite{MR3723158},
which established
the
existence of minimizers 
for~${{ {\hbox{{\Fontamici J}}} }}_{\Omega}$
in dimension~8 and above whose level
sets are asymptotic to a non-trivial cone (as a consequence,
the result
in~\eqref{89019230234855757757583} is optimal
in terms of dimension). Using a method introduced in~\cite{MR2140525},
the construction in~\cite{MR3723158} also provides, as a byproduct,
several examples 
in dimension~9 and higher which are structurally different from the one in~\cite{MR2846486}.
See also~\cite{MR2538506} for a detailed study of the Morse index
of the saddle-shaped solutions of the Allen-Cahn equation.

\section{The fractional Allen-Cahn equation}\label{S2}

Now, we consider a nonlocal analogue of the Allen-Cahn equation in~\eqref{AC1}
and we investigate which of the above classical results remain valid
also in this generality. The model that we take into account aims at
dealing with long-range interactions which can influence the coexistence
of these two phases, and indeed these contributions ``coming from
far'' can (and typically do) produce a number of new phenomena
with respect to the classical case.\medskip

The problem that we discuss involves the fractional Laplace operator
with~$s\in(0,1)$, defined as
$$ (-\Delta)^s u(x):=c_{n,s}\int_{\R^n}\frac{2u(x)-u(x+y)-u(x-y)}{|y|^{n+2s}}\,dy.$$
Here~$c_{n,s}$ is a suitable
renormalization constant, chosen in such a way that, if~$u$
is smooth and rapidly decreasing, then the Fourier Transform
of~$(-\Delta)^s u$ coincides with $|\xi|^{2s} \hat u(\xi)$, being~$\hat u$
the Fourier Transform of~$u$ (equivalently, the fractional Laplacian
acts as multiplication by~$|\xi|^{2s}$ in Fourier space).
See e.g.~\cites{MR0214795, MR0290095, MR2707618, MR2944369, MR3469920, 2017arXiv171011567A, 2017arXiv171203347G}
for the basics on this operator and for several
motivations and applications.
The fractional counterpart of the Allen-Cahn equation in~\eqref{AC1}
is
\begin{equation}\label{ACs}
(-\Delta)^s u = u-u^3\quad{\mbox{ in }}\Omega.
\end{equation}
Akin its classical counterpart, equation~\eqref{ACs} also comes
from a variational principle and it corresponds to the minimization
of the energy functional
\begin{equation}\label{Es}
{{ {\hbox{{\Fontamici J}}} }}_{s,\Omega}(u):=\frac{c_{n,s}}{4}\,\iint_{Q_\Omega}
\frac{\big| u(x)-u(y)\big|^2}{|x-y|^{n+2s}}\,dx\,dy+
\int_\Omega W\big(u (x)\big) \,dx,
\end{equation}
up to scaling constants which are omitted for the sake of simplicity.
In~\eqref{Es} we used the notation
\begin{eqnarray*}
Q_\Omega&:=& \R^{2n}\setminus (\R^n\setminus\Omega)^2\\&=&
\Big( \Omega\times\Omega\Big)\cup
\Big( \Omega\times(\R^n\setminus\Omega)\Big)\cup
\Big( (\R^n\setminus\Omega)\times\Omega\Big).
\end{eqnarray*}
Comparing the nonlocal energy in~\eqref{Es} with its classical counterpart
in~\eqref{E1}, we notice that the fractional model takes into account
long-range particle interactions. As a matter of fact, the local interfacial
term modeled in~\eqref{E1} by the Dirichlet energy is replaced
in~\eqref{Es} by a Gagliardo-Sobolev-Slobodeckij seminorm.
The role of the domain~$Q_\Omega$ is to collect all the couples~$(x,y)\in\R^n\times\R^n$
for which at least one of the points~$x$, $y$ belongs to
the container~$\Omega$:
interestingly, while in~\eqref{E1} the interface term takes into account
the points in~$\Omega$, which can be seen as the complement in~$\R^n$
of
the ``inactive'' set~$\R^n\setminus\Omega$,
the domain~$Q_\Omega$ describing the long-range interaction in~\eqref{Es}
consists in the complement in~$\R^{2n}$
of the ``inactive'' couples of points in~$(\R^n\setminus\Omega)^2$
(hence, $Q_\Omega$ takes into account all the ``active'' points
which interact with points inside the container).\medskip

In recent years, a great amount of research was carried out
on equation~\eqref{ACs} and on the energy functional~\eqref{Es} (and, in fact,
also other types of long-range interactions have been taken into account, see~\cite{MR1612250}).
Our goal here is to describe the results
obtained in this fractional framework which are related (and possibly similar
in spirit, or significantly different) with the ones described in Section~\ref{S1}.
\medskip

First of all, we point out that an analogue of the classical
$\Gamma$-convergence result holds true in the nonlocal setting,
with an important variant: namely, the limit in~\eqref{limit E}
remains valid, but the limit set~$E$ has different features, according
to the fractional parameter~$s$. 
We will state this result in the forthcoming Theorem~\ref{SVi134}.
To this end, to treat the case~$s\in(0,1/2)$, we need to recall
the notion of fractional minimal surface, as introduced in~\cite{MR2675483}.
Given~$s\in(0,1/2)$ and two (measurable) disjoint
sets~$A$, $B\subseteq\R^n$,
one defines the $s$-interaction between~$A$ and~$B$ by
$$ I_s(A,B):=\iint_{A\times B}\frac{dx\;dy}{|x-y|^{n+2s}}.$$
Then, the $s$-perimeter of a set~$E$ in the domain~$\Omega$
is defined as
$${\operatorname{Per}}_s(E,\Omega):=I_s(E\cap\Omega, E^c\cap\Omega)+
I_s(E\cap\Omega, E^c\cap\Omega^c)+I_s(E\cap\Omega^c,E^c\cap\Omega),
$$
where we used the complimentary set notation~$A^c:=\R^n\setminus A$.

To state the $\Gamma$-convergence result, it is also convenient
to introduce a scaled version of the fractional functional in~\eqref{Es}.
Namely, we set
$$ {{ {\hbox{{\Fontamici J}}} }}_{s,\Omega,\e}(u):
=\left\{
\begin{matrix}
\e^{2s-1}\displaystyle
\iint_{Q_\Omega}\displaystyle
\frac{\big| u(x)-u(y)\big|^2}{|x-y|^{n+2s}}\,dx\,dy+\displaystyle\frac1\e
\int_\Omega W\big( u(x)\big) \,dx &
{\mbox{ if }}s\in \left(\displaystyle\frac12,1\right),\\
\displaystyle\frac1{|\log\e|}\iint_{Q_\Omega}
\displaystyle\frac{\big| u(x)-u(y)\big|^2}{|x-y|^{n+2s}}\,dx\,dy+
\displaystyle\frac1{\e\,|\log\e|}\int_\Omega W\big( u(x)\big)\,dx&
{\mbox{ if }}s=\displaystyle\frac12,\\
\displaystyle\iint_{Q_\Omega}
\displaystyle\frac{\big| u(x)-u(y)\big|^2}{|x-y|^{n+2s}}\,dx\,dy+
\displaystyle\frac1{\e^{2s}}\int_\Omega W\big( u(x)\big)\,dx&
{\mbox{ if }}s\in \left(0,\displaystyle\frac12\right).
\end{matrix}
\right. $$
This scaled functional is obtained from~\eqref{Es}, using the blow-down
sequence in \eqref{BLSDOW}, after a multiplication that keeps
the energy of the one-dimensional
profile bounded uniformly in~$\e$.
In this setting, we can state the $\Gamma$-convergence
result for the fractional Allen-Cahn functional as follows:

\begin{theorem}[Theorem 1.5 in \cite{MR2948285}]\label{SVi134}
If~$\Omega$ is a smooth domain and~$u_\e:\Omega\to[-1,1]$ is 
a sequence of minimizers for~${{ {\hbox{{\Fontamici J}}} }}_{s,\Omega,\e}$
such that
$$\sup_{\e\in(0,1)}{{ {\hbox{{\Fontamici J}}} }}_{s,\Omega,\e}(u_\e)<+\infty,$$
then, up to a subsequence,
\begin{equation*} \lim_{\e\searrow0} u_\e=u_0:=\chi_E-\chi_{\R^n\setminus
E}\quad{\mbox{ in }}
L^1(\Omega),\end{equation*}
for some set~$E\subseteq\R^n$.

If~$s\in[1/2,1)$, the set~$E$
minimizes the perimeter in~$\overline\Omega$
with respect to its boundary datum.

If instead~$s\in(0,1/2)$
and~$u_\e$ converges weakly to~$u_0$ in~$\R^n\setminus\Omega$, then
the set~$E$ minimizes the fractional perimeter~${\operatorname{Per}}_s$
in~$\Omega$ with respect to its datum in~$\R^n\setminus\Omega$.
\end{theorem}

It is interesting to remark that the $\Gamma$-convergence results in Theorem~\ref{SVi134}
are significantly easier in the case~$s\in(0,1/2)$, since characteristic functions
are admissible competitors, having finite energy. Instead, the case~$s\in[1/2,1)$
is much harder to treat, since one has to ``reconstruct'' a local energy
from all the nonlocal contributions in the limit, and therefore a fine
measure theoretic analysis of integral contributions is needed in this case.

We now briefly discuss the fractional version of the
energy and density estimates.
We will see that
the estimate in~\eqref{EN:DE2}, which is somehow geometric
(stating that, in the large, both phases occupy a non-negligible portion
of a large ball centered at the interface), still holds in the fractional case.
Conversely, the energy bound in~\eqref{EN:DE1}
is influenced by the fractional parameter~$s$, in the same way
as the one presented in Theorem~\ref{SVi134}: indeed, for large values
of the parameter~$s$, the estimate in~\eqref{EN:DE1}
remains the same, while for small values of~$s$ the energy contributions
``coming from infinity'' are not anymore negligible and
they carry an additional amount of energy in a large ball (though this
energy produced by the phase transition
remains negligible with respect to the size of the ball).
The precise result goes as follows:

\begin{theorem}[Theorems 1.3 and 1.4 in \cite{MR3133422}]\label{101-098t89204}
If~$R\ge1$ and~$
u:B_{R+1}\to[-1,1]$ is a minimizer of~$
{{ {\hbox{{\Fontamici J}}} }}_{s,B_{R+1}}$,
then
\begin{equation*}
{{ {\hbox{{\Fontamici J}}} }}_{s,B_{R+1}}(u)\le \left\{\begin{matrix}
CR^{n-1} & {\mbox{ if }}s\in \left(\displaystyle\frac12,1\right),\\
CR^{n-1}\,\log R &{\mbox{ if }}s=\displaystyle\frac12,\\
CR^{n-2s} &{\mbox{ if }}s\in \left(0,\displaystyle\frac12\right),
\end{matrix}\right.\end{equation*}
for some~$C>0$.

In addition, if~$u(0)=0$, then
\begin{equation*}\begin{gathered}
{\mbox{the Lebesgue measures of $\{u>1/2\}$ and~$\{u<-1/2\}$
in~$B_R$}}\\{\mbox{are both greater than~$cR^n$,}}\end{gathered}
\end{equation*}
for some~$c>0$.
\end{theorem}

Of course, the constants in Theorem~\ref{101-098t89204} depend, in general,
on~$n$
and~$s$. Though weaker (at least for small~$s$)
than in the classical case, the estimates in Theorem~\ref{101-098t89204}
are sufficient to obtain the locally uniform convergence of the
level sets of minimizers,
as stated in the following result:

\begin{corollary}[Corollary 1.7 in \cite{MR3133422}]
If~$\Omega$ is a smooth domain, $E\subseteq\R^n$ and~$
u_\e:\Omega\to[-1,1]$
is a minimizer
of~${{ {\hbox{{\Fontamici J}}} }}_{s,\Omega,\e}$ such that~\eqref{limit E} holds true,
then the set~$\{ |u_\e|\le1/2\}$ converges locally uniformly in~$\Omega$
to~$\partial E$ as~$\e\searrow0$.\end{corollary}

Now, we discuss the fractional analogue of~\eqref{MODICA}.
To this aim, it is convenient to introduce the notion of
extension solution of~\eqref{ACs}
(see~\cite{MR2354493}). Namely, we consider the
Poisson Kernel
$$ \R^n\times(0,+\infty)=:\R^{n+1}_+\ni(x,t)\longmapsto
P(x,t):= \bar c_{n,s}\,\frac{t^{2s}}{(|x|^2+t^2)^{\frac{n+2s}2}},$$
where~$\bar c_{n,s}>0$ is the normalizing constant for which
$$ \int_{\R^n} P(x,t)\,dx=1,$$
for any~$t>0$. Given~$u:\R^n\to[-1,1]$, we define
$$ \R^{n+1}_+\ni(x,t)\longmapsto E_u(x,t):=\int_{\R^n} P(x-y,t)\,u(y)\,dy.$$
Then, if~$u$ is sufficiently smooth,
we have that~$E_u$ reconstructs the fractional Laplacian of~$u$ as a weighted
Neumann term: more precisely, one has that~$E_u$ satisfies
\begin{equation*}
\left\{
\begin{matrix}
{\operatorname{div}}(t^\alpha\nabla E_u)=0 & {\mbox{ in }}\R^{n+1}_+,\\
\tilde c_s\,\displaystyle\lim_{t\searrow0} t^\alpha\partial_t E_u=-(-\Delta)^s u & {\mbox{ in }}\R^n,
\end{matrix}
\right.
\end{equation*}
where~$\alpha:=1-2s\in(-1,1)$.
The constant~$\tilde c_s>0$ is needed just for normalization purposes
(and it can be explicitly calculated, see e.g.
Remark~3.11(a) in~\cite{MR3165278}).
Hence, if~$u$ is a solution of~\eqref{ACs},
then~$E_u$ is a solution of 
\begin{equation}\label{EXTE ACs}
\left\{
\begin{matrix}
{\operatorname{div}}(t^\alpha\nabla E_u)=0 & {\mbox{ in }}\R^{n+1}_+,\\
\tilde c_s\,\displaystyle\lim_{t\searrow0} t^\alpha\partial_t E_u=u^3-u & {\mbox{ in }}\R^n.
\end{matrix}
\right.
\end{equation}
Since, to the best of our knowledge, the
fractional counterpart of~\eqref{MODICA} is at the moment understood
only when~$n=1$, we will consider in~\eqref{EXTE ACs}
the case in which~$x\in\R$, namely
\begin{equation*}
\left\{
\begin{matrix}
{\operatorname{div}}(t^\alpha\nabla E_u)=0 & {\mbox{ in }}\R^{2}_+,\\
\tilde c_s\,\displaystyle\lim_{t\searrow0} t^\alpha\partial_t E_u=u^3-u
& {\mbox{ in }}\R,
\end{matrix}
\right.
\end{equation*}
and look at the related energy functional
$$ F(x,y):=(1-s)\,\int_0^y t^\alpha
\,\Big( |\partial_x E_u(x,t)|^2-|\partial_t E_u(x,t)|^2\Big)\,dt.$$
In this setting, the following result holds
true:

\begin{theorem}[Theorem~2.3(i) of~\cite{MR3165278}]\label{9Xuewihf383883}
Let~$u:\R\to[-1,1]$ be a solution of~\eqref{ACs}
such that
$$ \partial_{x} u(x)>0 \qquad{\mbox{and}}\qquad
\lim_{x\to\pm\infty} u(x)=\pm1.$$
Then, for any~$x\in\R$ and any~$y\ge0$ we have that
$$ F(x,y)\le W(u(x))=F(x,+\infty).$$
\end{theorem}

Interestingly, semilinear 
fractional equations possess
a formal Hamiltonian structure in infinite dimensions (see
Section~1.1 in~\cite{MR3165278}) and Theorem~\ref{9Xuewihf383883}
recovers
the classical Conservation of Energy Principle as~$s\nearrow1$
(see Section~6 in~\cite{MR3165278}).
It is an open problem to understand the possible validity of
results as in Theorem~\ref{9Xuewihf383883} when~$n\ge2$.\medskip

The last part of this note aims at discussing the recent developments of
the symmetry results for solutions of equation~\eqref{ACs}, in view
of the problems posed in
Conjectures~\ref{GDD} and~\ref{GDD2}
for the classical case.
As a matter of fact, the analogue of Conjecture~\ref{GDD2}
possesses a positive answer also in the fractional setting,
for any dimension~$n$ and any fractional exponent~$s\in(0,1)$,
see Theorem~2 in~\cite{MR2952412}.\medskip

As for the analogue of Conjecture~\ref{GDD} in the fractional framework, 
the problem is open in its generality, but it possesses a positive
answer for all~$n\le3$ and~$s\in(0,1)$, and also for~$n=4$ and~$s=1/2$,
according to the following result:

\begin{theorem}\label{898989}
Let~$u:\R^n\to[-1,1]$ be a solution of~\eqref{ACs}
in the whole of~$\R^n$ such that
\begin{equation*}
\frac{\partial u}{\partial x_n}(x)>0 \quad{\mbox{ for all }}x\in\R^n.
\end{equation*}
Suppose that either
$$ n\le 3 \qquad{\mbox{and}}\qquad s\in(0,1),$$
or
$$ n= 4 \qquad{\mbox{and}}\qquad s=\frac12.$$
Then~$u$ is $1$D.
\end{theorem}

Theorem~\ref{898989} is due to~\cite{MR2177165} when~$n=2$
and~$s=1/2$, \cites{MR2498561, MR3280032} when~$n=2$ and~$s\in(0,1)$,
\cite{MR2644786} when~$n=3$ and~$s=1/2$,
\cite{MR3148114} when~$n=3$ and~$s\in(1/2,1)$,
\cite{MR3740395} 
(based also on preliminary rigidity results in~\cite{2016arXiv161110105D})
when~$n=3$ and~$s\in(0,1/2)$,
\cite{2017arXiv170502781F} when~$n=4$ and~$s=1/2$.
The cases remained open will surely provide several very interesting
and challenging complications. It is also worth to point out that,
at the moment, there is no counterexample in the literature
to statements
as the one in Theorem~\ref{898989} in higher dimensions --
nevertheless
an important counterexample to the validity of Theorem~\ref{898989}
in dimension~$n\ge 9$ when~$s\in(1/2,1)$
has been recently announced by
H. Chan, J. D\'avila, M. del Pino, Y. Liu and J. Wei 
(see the comments after Theorem~1.3
in~\cite{2017arXiv171103215C}).\medskip

The validity of Theorem~\ref{898989} in higher dimensions
under the additional limit assumption in~\eqref{LIMINU}
has been also investigated in the recent literature.
At the moment, the best result known on this topic can be summarized as follows:

\begin{theorem}\label{9292019292:p202}
Let~$n\le8$.
Then, there exists~$\e_0\in\left(0, \frac12\right]$
such that for any~$ s\in\left(\frac12-\epsilon_0,1\right)$
the following statement holds
true.

Let~$u:\R^n\to[-1,1]$ be a solution of~\eqref{ACs}
in the whole of~$\R^n$ such that
\begin{equation*}\begin{gathered}
\frac{\partial u}{\partial x_n}(x)>0 \quad{\mbox{ for all }}x\in\R^n\\
{\mbox{and}}\quad
\lim_{x_n\to\pm\infty} u(x',x_n)=\pm1\quad{\mbox{ for all }}x'\in\R^{n-1}.
\end{gathered}
\end{equation*}
Then, $u$ is $1$D.
\end{theorem}

Theorem~\ref{9292019292:p202} consists in fact of the superposition
of three different results, also obtained with a different approach.
The result of Theorem~\ref{9292019292:p202} when~$s$
is larger than~$1/2$ follows from Theorem~1.3 in~\cite{2016arXiv161009295S}.
When~$s=1/2$, the result was announced after Theorem~1.1
in~\cite{2016arXiv161009295S} and established in Theorem~1.3
of~\cite{2018arXiv180201710S}.
The case~$s\in\left(\frac12-\epsilon_0,\frac12\right)$
has been established in Theorem~1.6 of~\cite{2016arXiv161110105D}.
In this latter framework, the quantity~$\epsilon_0$
is a universal constant (unfortunately, not
explicitly computed by the proof), and the arguments of the proof rely on
it in order to deduce the flatness of the corresponding limit interface,
which is in this case described by nonlocal minimal surfaces:
since such flatness results are only known above the threshold
provided by~$\epsilon_0$ (see Theorems~2--5
in~\cite{MR3107529}),
also Theorem~\ref{9292019292:p202} suffers of this restriction.
Of course, it is an important open problem to establish whether Theorem~\ref{9292019292:p202}
holds true for a wider range of fractional parameter,
as well as it would be very interesting to establish optimal regularity
results for nonlocal minimal surfaces.\medskip

The fractional counterpart of classical symmetry results
under the minimality assumption in~\eqref{89019230234855757757583}
has also been taken into account, with results
similar to Theorem~\ref{9292019292:p202}, which can be
summarized as follows:

\begin{theorem}\label{9eor7429942hd383737}
Let~$n\le7$.
Then, there exists~$\e_0\in\left(0, \frac12\right]$
such that for any~$ s\in\left(\frac12-\epsilon_0,1\right)$
the following statement holds
true.

Let~$u:\R^n\to[-1,1]$ be a minimizer
for~${{ {\hbox{{\Fontamici J}}} }}_{s,\Omega}$
in any bounded domain~$\Omega\subset\R^n$. Then, $u$ is $1$D.
\end{theorem}

Once again, Theorem~\ref{9eor7429942hd383737}
is a collage of different results obtained by different
methods and dealing with different parameter ranges.
Namely, the statement in Theorem~\ref{9eor7429942hd383737}
when~$s$ is larger than~$1/2$ has been proved in Theorem~1.2
of~\cite{2016arXiv161009295S}, and the case~$s=1/2$
has been treated in Theorem~1.2 of~\cite{2018arXiv180201710S}.
The case~$ s\in\left(\frac12-\epsilon_0,\frac12\right)$
has been established 
in Theorem~1.5 of~\cite{2016arXiv161110105D} (once again,
in this context, the threshold given by~$\epsilon_0$ is used
to apply the regularity results for nonlocal minimal surfaces in~\cite{MR3107529}
and it is a very interesting problem to determine
the possible validity of Theorem~\ref{9eor7429942hd383737}
when the dimensional and quantitative conditions are violated).\medskip

We think that it is important to stress the fact that the differences
between the fractional exponent ranges in the previous results
do not reflect a series of merely technical difficulties, but instead
it reveals fundamental structural differences between the phase
transitions when~$s\in[1/2,1)$ and when~$s\in(0,1/2)$.
These differences are somehow inherited by the dichotomy provided
in Theorem~\ref{SVi134}: indeed, as pointed out in this result,
when~$s\in[1/2,1)$ the nonlocal phase transitions end up showing
an interface corresponding to a local problem, while when~$s\in(0,1/2)$
the nonlocal features of the problem persist at any scale and produce
a limit interface of nonlocal nature. The structural differences
between local and nonlocal minimal surfaces may therefore
produce significant differences on the phase transitions too:
as a matter of fact, it happens that when~$s\in(0,1/2)$
the long-range interactions of points of the interface provide a number
of additional rigidity properties which have no counterpart
in the classical case. To exhibit a particular phenomenon related
to this feature, we recall the forthcoming result in Theorem~\ref{0202o33o31001}.
To state this result in a concise way, we introduce the notion
of ``asymptotically flat'' interface, which can be stated as follows.
First of all, we say that the interface of~$u$ in~$B_R$
is trapped in a slab of width~$2aR$ in direction~$\omega\in S^{n-1}$ if
\begin{equation}\label{AS FLAT}
\begin{split}
& \{ x\in B_R {\mbox{ s.t. }} \omega\cdot x\le -aR\}\subseteq
\left\{x\in B_R {\mbox{ s.t. }} u(x)\le -\frac{9}{10} \right\}
\\{\mbox{and }}\;&
\left\{x\in B_R {\mbox{ s.t. }} u(x)\le \frac{9}{10} \right\}\subseteq
\{ x\in B_R {\mbox{ s.t. }} \omega\cdot x\le aR\}.
\end{split}
\end{equation}
Of course, when~$a\ge1$, such condition is always satisfied, but the smaller the~$a$ is,
the flatter the interface is in the ball~$B_R$.
We say that the interface of~$u$ is asymptotically flat if there exists~$R_0>0$
such that for any~$R\ge R_0$
there exist~$\omega(R)\in S^{n-1}$ and~$a(R)\ge0$ such that
the interface of~$u$ in~$B_R$
is trapped in a slab of width~$2a(R)\,R$ in direction~$\omega(R)$ with
$$ \lim_{R\to+\infty} a(R)=0.$$
Roughly speaking, the interface of~$u$ is asymptotically flat
if, in large balls, it is trapped into slabs
with small ratio between the
width of the slab and the radius of the ball
(possibly, up to rotations which can vary from one scale to another).
In this setting, we have:

\begin{theorem}[Theorem 1.2 in~\cite{2016arXiv161110105D}]\label{0202o33o31001}
Let~$ s \in (0, 1/2)$
and~$u$ be a solution of~\eqref{ACs}
in~$\R^n$. Assume that the
interface of~$u$ is asymptotically flat.
Then, u is $1$D.
\end{theorem}

We think that Theorem~\ref{0202o33o31001} reveals several surprising aspects
of nonlocal phase transitions in the regime~$s\in(0,1/2)$, where
the contributions from infinity happen to be dominant. Indeed,
the result in Theorem~\ref{0202o33o31001} is valid for all solutions,
without any monotonicity or energy restrictions. This suggests that
if one has a phase coexistence in this regime, plugging additional energy into the system
can only produce two alternatives:
\begin{itemize}
\item either the interface oscillates significantly at infinity (i.e., the flatness
assumption
of Theorem~\ref{0202o33o31001} is not satisfied),
\item or the graph of the function~$u$ that
describes the state parameter of the system can
oscillate, but (due to Theorem~\ref{0202o33o31001}) such function is necessarily $1$D
and therefore the phase separation occurs along parallel hyperplanes, with possible multiplicity. 
\end{itemize}
It is also interesting to stress that a result as the one in Theorem~\ref{0202o33o31001}
does not hold for the classical Allen-Cahn equation (and indeed
Theorem~\ref{0202o33o31001} reveals a purely nonlocal phenomenon). As a matter of fact,
in Theorem~1 of~\cite{MR3019512} a solution of~\eqref{AC1} in~$\R^3$
is constructed 
whose level sets resemble an appropriate
dilation of a catenoid:
namely, the level sets of this solution lie in the 
asymptotically flat region~$\{ x=(x',x_3)\in\R^3
{\mbox{ s.t. }}|x_3|\le C(1+\log(1+|x'|)\}$, for a suitable~$C>0$. In particular,
condition~\eqref{AS FLAT} is satisfied with~$\omega(R):=(0,0,1)$
and~$a(R):=\frac{C(1+\log(1+R)}{R}$, which is infinitesimal as~$R\to+\infty$ and,
as a byproduct, 
the interface of this solution is asymptotically flat. Clearly, the solution constructed
in~\cite{MR3019512} is not $1$D, since its level sets are modeled on a catenoid rather
than on a plane, and therefore this example shows that an analogue of
Theorem~\ref{0202o33o31001} is false in the classical case.

A fractional counterpart of~\cite{MR3019512} has been recently provided
in~\cite{2017arXiv171103215C}, in the fractional regime~$s\in(1/2,1)$.
In particular,
Theorem~1.3 of~\cite{2017arXiv171103215C} establishes the existence of
an entire solution of~\eqref{ACs} in~$\R^3$
vanishing
on a rotationally symmetric surface which resembles a catenoid
with sublinear growth at infinity. This example shows that
an analogue of
Theorem~\ref{0202o33o31001} is false when~$s\in(1/2,1)$.

At the moment, it is an open problem to construct solutions of~\eqref{ACs} in~$\R^3$
with level sets modeled on a catenoid when~$s=1/2$, see Remark~1.4
in~\cite{2017arXiv171103215C}: on the one hand, the case~$s=1/2$ relates
the large-scale picture of the interfaces to the classical (and not to the nonlocal)
minimal surfaces (recall Theorem~\ref{SVi134}), therefore it is still conceptually
possible to construct catenoid-like examples in this setting; on the other hand,
the infinite dimensional gluing method in~\cite{2017arXiv171103215C} deeply relies
on the condition~$s\in(1/2,1)$, therefore important modifications would be needed
to achieve similar results when~$s=1/2$.

Interestingly, nonlocal catenoids corresponding to the case~$s\in(0,1/2)$
have been constructed in~\cite{2014arXiv1402.4173D} but, remarkably,
such surfaces possess linear (rather than sublinear)
growth at infinity (therefore, possible solutions of~\eqref{ACs}
modeled on such catenoids would not possess asymptotically flat interfaces,
which is indeed in agreement with Theorem~\ref{0202o33o31001}).
\medskip

We end this note with a few comments on the proof of Theorem~\ref{0202o33o31001}:
the main argument is an ``improvement of flatness''
which says that if a sufficiently sharp interface
is appropriately flat ``from the unit ball~$B_{1}$ towards infinity'', then
it is even flatter in~$B_{1/2}$ (see Theorem~1.1
in~\cite{2016arXiv161110105D} for full details).
Suitable iterations of this argument
give a control of the interface all the way to infinity, showing in particular
that (possibly after a rotation)
the interface is trapped between a graph that is
Lipschitz and sublinear at infinity and its translate.
This control of the growth at infinity of the interface in turn allows
the use of the sliding method ``in a tilted direction''. Namely,
one fixes~$e'\in\R^{n-1}$ with~$|e'| = 1$ and~$\delta>0$ and set
$$ e_\delta:=\frac{(e',\delta)}{\sqrt{1+\delta^2}}\in S^{n-1}\qquad{\mbox{and}}\qquad
u^{(t)}(x):= u(x-e_\delta t).$$
We point out that~$u^{(t)}$ is the translation of the original solution~$u$
in the slightly oblique direction~$e_\delta$ and so the growth control of the
interface, combined with a precise estimate of the decay of the solution
and the maximum principle,
implies that~$u^{(t)}$ lies below~$u$ for $t$ sufficiently large (say, $t\ge T(e',\delta)$,
and we observe that the use of maximum principle here
relies on the monotonicity
property of the Allen-Cahn nonlinearity outside the interface, namely the function~$f(r):=r-r^3$
is decreasing when~$|r-1|\le\frac1{10}$).

Then, one keeps sliding~$u^{(t)}$, reducing the value of~$t$, and using again the
maximum principle it follows that~$u^{(t)}\le u$ for any~$t\ge0$.
As a consequence of this, for any~$t\ge0$,
any~$x=(x',x_n)\in\R^n$ and any~$e'\in S^{n-2}$,
$$ 
u\left((x',x_n)-\frac{(e't,\delta t)}{\sqrt{1+\delta^2}}\right)=
u(x-e_\delta t)=
u^{t}(x)\le u(x)$$
and accordingly, sending~$\delta\searrow0$,
\begin{equation}\label{90x1q0}
u(x'-e' t,x_n)\le u(x) .
\end{equation}
Writing~\eqref{90x1q0} with~$e'$ replaced by~$-e'$ (as well as~$x$ replaced by~$y$), it follows that
\begin{equation}\label{90x1q1}
u(y'+e' t,y_n)\le u(y) ,
\end{equation}
for any~$y\in\R^n$ and any~$e'\in S^{n-2}$.
Then, choosing~$y:=x-(e't,0)$ in~\eqref{90x1q1} and using again~\eqref{90x1q0},
$$ u(x)=
u(x'-e't+e' t,x_n)\le u(x'-e't,x_n)\le u(x)$$
and therefore
$$ u(x)=u(x'-e't,x_n),$$
for every~$x\in\R^n$, every~$t\ge0$ and every~$e'\in S^{n-2}$. This shows that, possibly
after a rotation, the solution~$u$
depends only on~$x_n$ and so
it completes the proof of Theorem~\ref{0202o33o31001}.

\vfill

{\footnotesize

\noindent {\em Addresses:} \\

Serena Dipierro \& Enrico Valdinoci.
Dipartimento di Matematica, Universit\`a di Milano,
Via Saldini 50, 20133 Milan, Italy, and
School of Mathematics
and Statistics,
University of Western Australia,
35 Stirling Hwy, Crawley WA 6009, Australia.\\

Enrico Valdinoci.
School of Mathematics
and Statistics,
University of Melbourne, 813 Swanston St,
Parkville VIC 3010, Australia, and
IMATI-CNR, Via Ferrata 1, 27100 Pavia,
Italy. \\
\smallskip

\noindent{\em Emails:}\\

{serena.dipierro@unimi.it}, {enrico@mat.uniroma3.it}  

}

\end{document}